\documentclass{amsart}
\usepackage[T1]{fontenc}
\usepackage{tgtermes}
\usepackage[french,english,american]{babel}
\usepackage{amsmath,amssymb,amsthm}
\usepackage{mathrsfs}
\usepackage{graphicx}
\usepackage[all]{xy}
\usepackage{tikz,xcolor,color}
\usetikzlibrary{shapes.geometric}
\usepackage{setspace}
\newtheorem{theorem}{Theorem}[section]

\theoremstyle{definition}
\newtheorem{definition}{Definition}[section]

\newtheorem*{remark}{Remark}

\newtheorem*{thm*}{Theorem}
\newcommand{\set}[1]{\left\{#1\right\}}
\newcommand{\paren}[1]{\left(#1\right)}

\newcommand{\R}{\mathbb{R}}
\newcommand{\eng}{\mathscr{E}}
\newcommand{\eff}{R_{\operatorname{eff}}}
\newcommand{\alb}{\mathsf{X}}

\newcommand{\crits}{\mathcal{C}}
\newcommand{\pcrits}{\mathcal{P}}

\newcommand{\sss}{\mathcal{L}}

\newcommand{\shift}{\sigma}

\newcommand{\shiftspace}{\Sigma}
\newcommand{\leftshift}{\sigma}
\newcommand{\trace}{\operatorname{Tr}}

\def\pcf{{\mbox{p.c.f.}\xspace}}

\def\s-s{self-similar}
\numberwithin{equation}{section}
\title[Resistance scaling factor of fractalina and pillow fractals]{Resistance scaling factor of the pillow and fractalina fractals}
\author[M.~J.~Ignatowich]{Michael J.~Ignatowich}
\address[M.~J.~Ignatowich]{Department of Mathematics, University of Connecticut, Storrs, CT 06269, USA}
\email{Michael.Ignatowich@uconn.edu}
\author[D.~J.~Kelleher]{Daniel J.~Kelleher}
\address[D.~J.~Kelleher]{Department of Mathematics, University of Connecticut, Storrs, CT 06269, USA}
\email{daniel.kelleher@uconn.edu}
\author[C.~E.~Maloney]{Catherine E. Maloney}
\address[C.~E.~Maloney]{Department of Mathematics, University of Connecticut, Storrs, CT 06269, USA}
\email{Catherine.E.Maloney@uconn.edu}
\author[D.~J.~Miller]{David J. Miller}
\address[D.~J.~Miller]{Department of Mathematics, University of Connecticut, Storrs, CT 06269, USA}
\email{David.J.Miller@uconn.edu}
\author[K.~Nechyporenko]{Khrystyna Nechyporenko}
\address[K.~Nechyporenko]{Department of Mathematics, University of Connecticut, Storrs, CT 06269, USA}
\email{Khrystyna.Nechyporenko@uconn.edu}
\thanks{Research supported in part by NSF grant  DMS 0505622}
\begin{document}

\begin{abstract}
Much is known in the analysis of a finitely ramified self-similar fractal when the fractal has a harmonic structure: a Dirichlet form which respects the self-similarity of a fractal. What is still an open question is when such structure exists in general. In this paper, we introduce two fractals, the fractalina and the pillow, and compute their resistance scaling factor. This is the factor which dictates how the Dirichlet form scales with the self-similarity of the fractal. By knowing this factor one can compute the harmonic structure on the fractal. The fractalina has scaling factor $(3+\sqrt{41})/16$, and the pillow fractal has scaling factor $\sqrt[3]{2}$. 
\\[12pt]
\end{abstract}

\maketitle

\section{Introduction}
The first step in understanding analysis and probability on a fractal space is calculating a Dirichlet form on a suitable domain of functions. As seen in \cite{FOT11}, a Dirichlet form is equivalent to a Laplacian operator, and to a diffusion process on the fractal. The standard references for the theory of analysis and probability on fractals are \cite{Kig01, Str06, BN98}. However, techniques used in computing Dirichlet forms are often algebraic --- the current work rewrites the problem in terms of solving equations of homogeneous rational functions. This is the first step in using the theory of curves and surfaces to classify the space of forms on a given fractal.
 
We calculate the Dirichlet forms on two \s-s fractals, the modified fractalina Sierpinski gasket and the pillow fractal. Recently there has been much progress in establishing the existence and uniqueness of these forms, for example \cite{BSU08,CS07,BCF+07,HMT06,Pei08}. Despite this, there are still simple examples, such as the fractalina, which elude understanding. In \cite{NT08}, the theory of self-similar groups is used to establish the existence and uniqueness of a scaling factor for the pillow fractal. The current work follows up by determining the exact value of this factor.

Self-similarity --- the quality of looking the same on all scales --- is one of the most important concepts in the study fractal geometry and analysis on fractal spaces. A compact metric space $K$ is called self-similar if there exists a finite collection of injective contraction maps $\set{F_i}_{i=1}^N$ with the property that
\[
K = \bigcup_{i=1}^N F_i(K)
\]
Analysis on these fractals relies on Dirichlet forms. The classical example of a Dirichlet form is the energy integral $\int \nabla u \cdot\nabla v $ on $\Bbb{R}^d$. Dirichlet forms are bilinear forms the existence which is equivalent to Laplacian operators and diffusion processes. For more on this relation, see \cite{FOT11}. A Dirichlet form $\eng$ on a self-similar fractal is self-similar if for some $\rho_i>0$, $i=1,\ldots, N$
\begin{align}
\label{eq:ssform}\eng(u,u) = \sum_{i=1}^N \rho_i\eng(u\circ F_i,u\circ F_i)\ \text{.}
\end{align}
In the case where all of the $\rho_i=\rho$ are equal, we call $\rho^{-1}$ a resistance scaling factor. Calculating this factor allows us to much better understand the energy form on the limiting fractal. We calculate this factor for two specific self-similar structures, the modified fractalina Sierpinski gasket, or fractalina for short, and the pillow fractal.

One particularly nice property of these fractals is that they are post-critically finite (p.c.f.) self-similar structures --- there is a finite set $V_0$, called the fractal boundary, such that $F_i(K)\cap F_j(K) \subset F_i(V_0)\cap F_j(V_0)$, for all $i\neq j$. A theory of analysis exists for a \pcf\ \s-s structure with a harmonic structure, see \cite{Kig93,Kig01,Str06}. This literature also details the existence of a self-similar Dirichlet form for highly symmetric nested fractals.

The exact conditions for existence and uniqueness of a harmonic structure on self-similar fractals is still not fully understood, even for p.c.f. self-similar structures. There has been much progress in recent years, see \cite{Met03,Sab97}. In \cite{HMT06}, the symmetry conditions are greatly relaxed. The existence of harmonic structures on other large classes of fractals is found in \cite{Pei08}, including generalized Vicsek sets.

The works \cite{CS07, BSU08} identify a class of energy forms on non-equilateral Sierpinski gaskets. These energy forms are indexed by a 2 dimensional smooth manifold. This is proven by showing that the self-similar scaling factors of the forms satisfy certain polynomial relations and are thus contained in an algebraic variety. Even though the energy forms are not unique, algebraic techniques were used to classify them.

Fractalina is a modified version of the Sierpinski gasket. Like the gasket, fractalina is generated by three contraction mappings, but one of the contraction features a $180^\circ$ rotation.  This rotation introduces a few complications. The boundary of this fractal now consists of 4 points, and the fixed point of the rotated contraction is no longer in the boundary of the fractal. In section \ref{sec:fractalina} we introduce the fractalina and compute its resistance scaling factor to be $(3+\sqrt{41})/16$. 

 In section \ref{sec:pillow} we introduce the pillow and calculate the scaling factor $\sqrt[3]{2}$. The solution $\sqrt[3]{2}$ is in interesting result because in most of the known examples, the resistance scaling factor is the solution to a quadratic equation. Establishing a scaling factor also proves that there is a self-similar resistance form on both of these fractals. 
 
 In section \ref{sec:ssstructures} we introduce rigorous definitions of \s-s structures, and discuss how to approximate energy forms on these structures by looking at finite subsets. In section \ref{sec:matrices} we reduce the problem of defining a self-similar Dirichlet form to that of solving a matrix equation. Thus  the problem becomes a matter of solving a finite number of rational equations, a problem which is accessible to computer algebra systems. Section \ref{sec:electric} discusses the techniques from electrical circuits which we use in the paper.
\\ \noindent
\textbf{ Acknowledgements: } The authors would like to thank Alexander Teplyaev for his insight and guidance.
%
\section{\s-s structures}\label{sec:ssstructures}
This section describes that basics of the theory of \s-s structures and their Dirichlet forms. For a more extensive treatment, see \cite{Kig01}.

Consider the set  $\alb = \set{1,\ldots, N}$, which we shall call the alphabet. Then $\alb^n = \set{w_1w_2\cdots w_n~|~w_i\in\alb}$
will be the words of length $n$. We take $\alb^* = \bigcup_{n=1}^\infty \alb^n $ to be the collection of all finite words, and $\shiftspace = \prod_{i=1}^\infty \alb$ to be the set of all right-infinite words, called the shift-space of $\alb$.

Give $\alb$ the discrete topology. Thus $\shiftspace$ inherits the product topology (i.e. the topology whose basis is sets of the form $\prod_{i=1}^\infty A_i$, such that there is some $M$ with $A_i = \alb$ for $i \geq M$). We define the maps $\leftshift_i:\shiftspace\to\shiftspace$ by $\leftshift_i(w) = iw$, for all $w\in\shiftspace$. We also define the map $\shift:\shiftspace\to\shiftspace$ by $\shift(w_1w_2w_3\ldots) = w_2w_3\ldots$, for $w_1w_2w_3\ldots\in\shiftspace$.

\begin{definition}\label{defn:sss}
Let $K$ be a compact metrizable space, let $\alb$ be a finite set for $F_i:K\to K$ continuous injections such that $K = \cup_{i\in\alb}F_i(K)$. We call the triplet $\sss=(K,\alb,\set{F_i}_{i\in\alb})$  a self-similar structure on $K$ if there is a continuous surjection $\pi:\shiftspace\to K$ such that $F_i\circ \pi = \pi\circ \leftshift_i$, holds for all $w\in\alb^*$.
\end{definition}

If $w = w_1w_2w_3\ldots w_r$, then we define $F_w := F_{w_1}\circ F_{w_2}\circ\cdots\circ F_{w_r}$ and $K_w = F_w(K)$.   

\begin{definition}
Let $\sss=(K, \alb, \set{F_i}_{i\in \alb})$ be a \s-s structure on $K$. The critical set of $\sss$ is defined 
\[
\crits(\sss) = \linebreak \pi^{-1} \paren{\bigcup_{i, j\in \alb, i \neq j} (K_i \cap K_j) } \ \text{,}
\] 
and the post-critical set 
\[
\pcrits(\sss) = \bigcup_{n\geq1} \sigma^n (C)\ \text{,}
\]
where $\sigma$ is the shift operator on $\alb^{-\omega}$. We abbreviate the critical and post-critical sets by $\crits$ and $\pcrits$ when there is no danger of confusion. A \s-s structure is called post critically finite (\pcf) if $\pcrits$ is a finite set.

We shall define the fractal boundary to be $V_0 := \pi(\pcrits)$, and recursively define $V_n = \cup_{w\in\alb^n} F_w(V_0)$, in particular $K_u\cap K_v \subseteq V_n$ for any $u\neq v \in \alb^n$.
\end{definition}

For a finite set $V$, define the vector space $\ell(V) := \set{u:V\to\R}$ with pointwise addition and scalar multiplication. $\ell(V)$ is the set of vectors indexed by $V$.
  
We use equation \ref{eq:ssform} to build a Dirichlet form on $K$ from forms on $V_0$. We define the conductance matrix $C$ whose rows and columns are indexed by $V_0$, to be non-negative symmetric with diagonal entries which are identically $0$. That is $C_{pq}\geq 0$, $C_{pq} = C_{qp}$, and $C_{pp}=0$ for all $p$ and $q\in V_0$.  We define a bilinear form $\eng_0:\ell(V_0)\times \ell(V_0) \to \R$ by
\[
\eng_0(u,v) = \frac12\sum_{p,q\in V_0} (u(p)-u(q))(v(p)-v(q))C_{pq}
\]
For and bilinear form, we define the induced quadratic form $\eng_0(u) := \eng_0(u,u)$.

For $\rho = (\rho_0,\ldots,\rho_{N-1})$, inductively define an energy form on $V_n$ by a map $\Psi_\rho$ which takes an energy from on $V_{n-1}$ and produces an energy form on $V_n$ by 
\[
\eng_n(u,v) := \Psi_\rho(\eng_{n-1})(u,v) : = \sum_{i=1}^{N} \rho_i\eng_{n-1}(u\circ F_i,v\circ F_i).
\]
We define the energy form $\eng$ on $K$
\[
\eng(u,v) = \lim_{n\to\infty}\eng_n(u|_{V_n},v|_{V_n})\text{,}
\]
and refer to set of $u$ such that $\eng(u)$ is well defined as the domain of $\eng$.

However, we can not chose $\rho$ and $C$ arbitrarily, we require some consistency between approximating energies $\eng_n$.  Recall $V_{n-1} \subset V_n$. Define the operator $\trace$ which takes a Diriclet form on $V_n$ and produces an energy form on $V_{n-1}$. For $u\in\ell(V_{n-1})$,
\[
\trace\eng_n(u,u) := \inf\set{\eng_n(g)~|~g\in\ell(V_n),~g|_{V_{n-1}}=u}
\]
The consistency condition we need is that $\trace\eng_1 = \eng_0$. If this is the case, we can find an element $\tilde{u}\in \ell(V_n)$ with $\tilde{u}|_{V_{n-1}} = u$ and $\eng_n(\tilde{u},\tilde{u}) = \trace\eng_{n-1}(u,u) = \eng_{n-1}(u,u)$. $\tilde{u}$ is called a harmonic extension of $u$, we shall discuss such extensions more in section \ref{sec:pillow}.

Determining a harmonic structure on a self-similar structure $K$ is equivalent to finding a connected matrix $C$ (We  call an $n\times n$ matrix $M$ connected if for all $i,j = 1,2,\ldots,n$, $i\neq j$, there is $i=i_1,i_2,\ldots,i_k=j$ such that $M_{i_\ell i_{\ell+1}}\neq 0$ for $\ell = 1,2,\ldots,k-1$.) and $\rho$ such that $\trace\circ\Psi_\rho(\eng_0) = \eng_0$.
%
\section{Energy Forms as Matrices}\label{sec:matrices}

For $n_0 < n_1$, we think of $\R^{n_0}\otimes\R^{n_0}$ as $n_0\times n_0$ matrices, $\R^{n_1}\otimes\R^{n_1}$ as $n_1\times n_1$ matrices. Consider the collection of $n_1\times n_0$-matrices $\set{\psi_i}_{i=1}^N$ with entries which are all zero except for exactly one 1 in each row and at most one 1 in each column. Thus $\set{\psi_i}_{i=1}^N$ are linear injections $\psi_i:\R^{n_0}\to\R^{n_1}$. Assume that
\[
\operatorname{span} \set{\psi_i(\R^{n_0})~|~i=1,2,\ldots,N} = \R^{n_1}
\]

Define the map
\[
\Psi:(0,\infty)^N\times\R^{n_0}\otimes \R^{n_0} \to \R^{n_1}\otimes \R^{n_1}
\]  
 by 
\[
\Psi(\rho,A) = \sum_{i=1}^N \rho_i \psi_i A \psi_i^T
.\]
We can also think of $\Psi$ as defined on simple tensors by
\[
\Psi(\rho,x\otimes y) = \sum_{i=1}^N \rho_i(\psi_i(x)\otimes \psi_i(y)),
\]
and extended linearly.

Elements of $\R^{n_1}\otimes \R^{n_1}$ can be written as block matrices
\[
M = \left(\begin{array}{cc}A & B \\C & D\end{array}\right)
\]
where $A\in\R^{n_0}\otimes \R^{n_0}$, $B\in \R^{n_0}\otimes \R^{n_1}$, $C\in\R^{n_1}\otimes\R^{n_0}$, and $D\in\R^{n_1-n_0}\otimes\R^{n_1-n_0}$.
We define the map 
\[
\trace:\R^{n_1}\otimes\R^{n_1}\to\R^{n_0}\otimes\R^{n_0}
\]
by  $$\trace(M) = A-BD^{-1}C.$$ This is called the Schur complement in general.
Finally, we define the map 
\[
\Lambda = \trace\circ \Psi:  (0,\infty)^N\times\R^{n_0}\otimes \R^{n_0}\to \R^{n_0}\otimes \R^{n_0}.
\]

In the fractal theory, $\ell(V_0)$ is $\R^{n_0}$ and $\ell(V_1)$ is $\R^{n_1}$, where $n_0 = |V_0|$ and $n_1=|V_1|$. $\psi_i$ is the map which extends a function in $\ell(F_i(V_0))$ to a function $\ell(V_1)$ with values identically 0 on $V_1\setminus F_i(V_0)$. There is a natural correspondence between bilinear forms and matrices given by mapping a matrix $M$ to a form
\[
\eng_M(x,y) = (Mx)\cdot y = y^TMx.
\]
The energy forms in the previous section correspond to matrices $M$ that have non-negative off-diagonal entries and satisfy $M{\bf1}=0$ where $\bf1$ is a vector with all entries equal to $1$. In particular, $C_{pq}$ are the off diagonal entries.

It is simple to see that $\Psi_\rho\eng_M = \eng_{\Psi(\rho,M)}$. It is shown in \cite{Kig01} chapter 2 that $\trace \eng_M = \eng_{\trace M}$, where $\trace \eng_M$ is defined as in the last section.

Thus finding a self-similar energy form is the same as finding $\rho\in(0,\infty)^N$ and a connected matrix $M\in \R^{n_0}\otimes\R^{n_0}$, as above, such that
\begin{align}
\Lambda(\rho,M) =M \ ?\label{eq:eigprob}
\end{align}

%
\section{Some basics of electric circuits}\label{sec:electric}
The techniques used in computing the harmonic structure on the fractalina gasket come from the basic theory of electrical networks. For a full introduction to these networks and their relationship with Dirichlet forms see \cite{DS84}.

Consider $V$, a finite collection of nodes and  a symmetric matrix of conductances  $C = \set{C_{pq}}_{p,q\in V}$ with positive entries. We often draw $V$ as a graph, where there is an edge between $p$ and $q$ if $C_{pq}>0$. Likewise if $p$ is not connected to $q$ in this graph, then $C_{pq} = 0$.

Think of $C_{pq}$ has being the strength of an electrical connection between $p$ and $q$. On the other hand, one can consider the resistance between two points: $R_{pq} = 1/C_{pq}$  if $C_{pq}>0$ and $R_{pq} = \infty$ if $C_{pq} = 0$. If we take a function $u: V \to \mathbb{R}$ then we can take the energy of the state 
\[
\eng(u)=\frac12\sum_{p,q\in V} (u(p)-u(q))^2C_{pq} \ \text{,}
\] 
We are most often interested in computing the effective resistance between points $p$ and $q$ in $V$. Informally, this is the resistance if one node was hooked to the positive end of a battery, while the other was hooked to the negative. Formally, is defined by
\[
\frac{1}{\eff (p,q)} = \inf\set{\eng(u,u)~|~u:V\to\R,~f(u)=1,~f(u)=0},
\]
for $p\neq q$ and $\eff(p,p) = 0$. It can be shown that $\eff$ is a metric on $V$ (see \cite{Kig01}). 

When $\eng$ is a self-similar resistance form, effective resistance scales with the self-similarity of the fractal. Thus if we can calculate the effective resistance between two points in the fractal on two levels of levels approximations, then we can obtain the resistance scaling factor. This is the technique we employ to calculate the scaling factor of the fractalina, and is the same strategy used in \cite{BCF+07} to calculate the resistance of fractal $N$-gaskets.

This reduces the problem of finding the resistance scaling factor to that of finding the resistance of an electrical network. We shall need two tools from the basic theory of electrical circuits to reduce the conductance networks --- Kirchhoff's laws and the $\Delta$-Y transform.

\begin{figure}[t]
\begin{tikzpicture}
\draw (90:2cm)--node[left]{$R_{ab}$} (210:2cm)node[below]{$b$} --node[below]{$R_{bc}$}(330:2cm)node[below]{$c$}--node[right]{$R_{ac}$}(90:2cm)node[above]{$a$};

\filldraw(90:2cm)  circle (2pt);
\filldraw(210:2cm)  circle (2pt);
\filldraw(330:2cm)  circle (2pt);
\end{tikzpicture}\hspace{.5in}
\begin{tikzpicture}
\draw (90:2cm)node[above]{$a$}--node[left]{$r_{a}$} (0,0) (210:2cm)node[below]{$b$} --node[above]{$r_{b}$}(0,0) (330:2cm)node[below]{$c$}--node[above]{$r_{c}$}(0,0);

\filldraw(90:2cm)  circle (2pt);
\filldraw(210:2cm)  circle (2pt);
\filldraw(330:2cm)  circle (2pt);
\filldraw (0,0)  circle (2pt);
\end{tikzpicture}
\caption{The $\Delta$-Y transform.}\label{fig:deltay}
\end{figure}
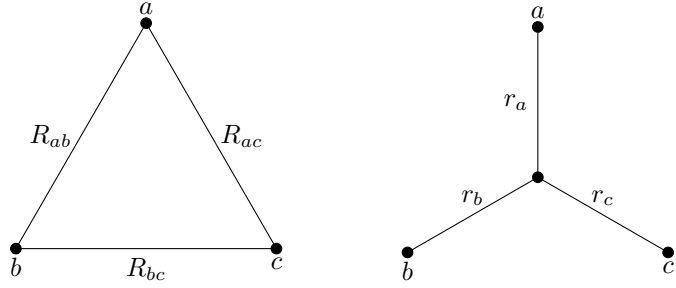
%
Kirchoff's two laws state that connections in parallel add conductances and connections in sequence add resistance. This allows us to reduce the graph as follows, if there are two connections between nodes $p$ and $q$ with conductances $C_{pq}$ and $C_{pq}'$, then the total conductance is $C_{pq} + C_{pq}'$. Likewise if there is a connection between node $p$ and $r$ with resistance $R_{pr}$ and between nodes $r$ and $q$ with resistance $R_{rq}$, and there are no other connections to $r$ (i.e. $C_{rs}=0$ for $s\in V$ with $s\neq p$ and $s\neq q$), then we can combine the two connections with a connection between $p$ and $q$ with resistance $R_{pq} = R_{pr}+R_{rq}$.

 If we can reduce a graph to one connection between $p$ and $q$, then the resistance of this connection is the effective resistance between these points. In particular, if we reduce the graph to be such that there is only one path between $p$ and $q$ which does not cross itself, then the effective resistance between $p$ and $q$ is the sum of the resistances of each edge of this path.

The other main tool we use in reducing networks is the $\Delta$-Y transformation. This tells us the process by which to transform a triangle network to a network which looks like a ``Y'' in a way that the effective resistance between any two nodes remains the same. If the resistances are set up as they are in figure \ref{fig:deltay}, then
\[
r_a = \frac{R_{ab}R_{ac}}{R_{ab}+R_{bc} + R_{ac}} \quad\text{and} \quad R_{ab} = \frac{r_ar_b + r_br_c+r_ar_c}{r_c}.
\]
This transformation is useful because it is easier to find the effective resistance between two nodes of a Y network than a $\Delta$ network. For example, $\eff(a,b) = r_a+r_b$.
%
%
\section{Fractalina Fractal} \label{sec:fractalina}
%
Fractalina was introduced by Robert Devaney using standard techniques in iterated function systems \cite{CDF98}. The construction of the fractalina fractal is similar to that of the Sierpinski triangle. 

Let $P_1,P_2$, and $P_3$ be vertices of an equilateral triangle in the complex plain, fractalina is the unique compact Hausdorff set invariant under the iterated function system (IFS) consisting of
\[
F_1(z) = \frac{z+P_1}{2},\quad F_2(z) = \frac{z+P_2}{2},\quad\text{and}\quad F_3=\frac{3P_3-z}{2}.
\]  
Similar to the IFS of the Sierpinski gasket the contractions $F_i$ have contraction factor of 1/2, however,  the contraction $F_3$ includes a rotation around $P_3$ of 180 degrees.

If assumptions are weaken so that $P_i$ are the vertices of any triangle, then the maps above result in a homeomorphic fractal. We assume equilateral, to allow us to make deductions based on the symmetry of the resulting fractal.
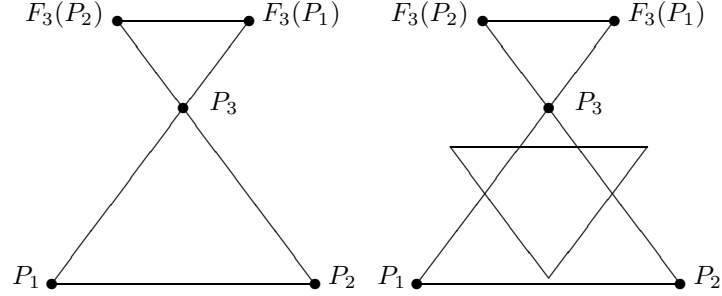
\begin{figure}[t]
\setlength{\unitlength}{3.5cm}
\begin{picture}(1,1)
\put(0,0){\line(3,4){.75}}
\put(-.15,0){$P_1$}
\put(.8,1){$F_3(P_1)$}
\put(-.1,1){$F_3(P_2)$}
\put(.6,.66){$P_3$}
\put(1,0){\line(-3,4){.75}}
\put(0,0){\line(1,0){1}}
\put(1.05,0){$P_2$}
\put(.25,1){\line(1,0){.5}}
\put(0,0){\circle*{.04}}
\put(1,0){\circle*{.04}}
\put(.25,1){\circle*{.04}}
\put(.75,1){\circle*{.04}}
\put(.5,.667){\circle*{.04}}
\end{picture}\hspace{.5in} 
\begin{picture}(1,1)
\put(0,0){\line(3,4){.75}}
\put(-.15,0){$P_1$}
\put(.8,1){$F_3(P_1)$}
\put(-.1,1){$F_3(P_2)$}
\put(.6,.66){$P_3$}
\put(1,0){\line(-3,4){.75}}
\put(0,0){\line(1,0){1}}
\put(.125,.52){\line(1,0){.75}}
\put(.125,.52){\line(3,-4){.375}}
\put(.875,.52){\line(-3,-4){.375}}
\put(1.05,0){$P_2$}
\put(.25,1){\line(1,0){.5}}
\put(0,0){\circle*{.04}}
\put(1,0){\circle*{.04}}
\put(.25,1){\circle*{.04}}
\put(.75,1){\circle*{.04}}
\put(.5,.667){\circle*{.04}}
\end{picture}

\caption{Level 1 and Level 2 of fractalina}\label{fig:linagraph}
\end{figure}

We assign resistances between each vertex of the Level 1 fractal, the graph hour-glass graph on the left in figure \ref{fig:linagraph}. In addition to the fractal boundary, the set $\set{P_1,P_2,F_3(P_1),F_3(P_2)}$, we include $P_3$ in our first level network, because removing this point disconnects the fractal. That is we assume
\[
V_0=\set{P_1,P_2,P_3,F_3(P_1),F_3(P_2)}\text{.}
\]
In addition, we assume the network consists of two triangular networks, $\set{P_1,P_2,P_3}$ and $\set{F_3(P_1),F_3(P_2),P_3}$. The second level is attained by applying the contractions, and we end up with the graph on the right in figure \ref{fig:linagraph}.

\begin{theorem}
The resistance scaling factor between any two levels of the fractalina is $(3+\sqrt{41})/16$.
\end{theorem}

\begin{figure}[t]
\setlength{\unitlength}{5cm}
\begin{picture}(1,1)
\put(0,0){\line(1,1){.75}}
\put(1,0){\line(-1,1){.75}}
\put(0,0){\line(1,0){1}}
\put(.25,.75){\line(1,0){.5}}
\put(.5, .78){$kR_1$}
\put(.7,.3){$R_2$}
\put(.6, .55){$kR_2$}
\put(.5,.03){$R_1$}
\end{picture}\hspace{.25in} 
\setlength{\unitlength}{5cm}
\begin{picture}(1,1)
\put(.5,.2){\line(0,1){.5}}
\put(.5,.2){\line(-3,-2){.2}}
\put(.5,.2){\line(3,-2){.2}}
\put(.5,.7){\line(-3,2){.1}}
\put(.5,.7){\line(3,2){.1}}
\put(.5, .78){$\alpha$}
\put(.52, .4){$\beta$}
\put(.55,.035){$\gamma$}
\end{picture}\hspace{.5in}
\caption{Delta-Y transformation applied to Level 1}\label{fig:level1deltay}
\end{figure}
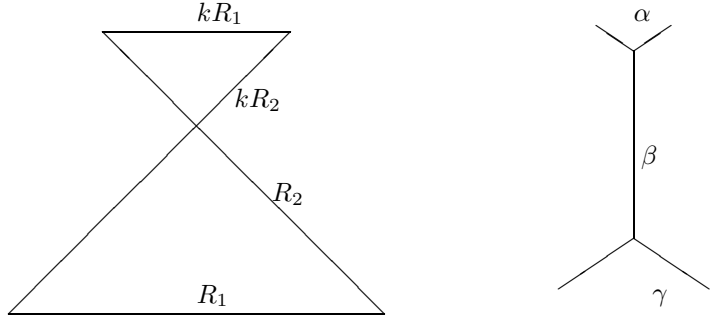

Let $k$ represent the scaling factor associated with the Fractalina contractions. Since
\[
F_3(\set{P_1,P_2,P_3}) = \set{F_3(P_1),F_3(P_2),P_3}\text{,}
\]
the resistances between elements on the right hand side are are those of $\set{P_1,P_2,P_3}$ scaled by $k$. Noting the reflection symmetry through the vertical axis, we arrive at the configuration of resistances found in figure \ref{fig:level1deltay}. This figure also shows the application of a $\Delta$-Y transformation to the two triangles. The three new resistances in figure \ref{fig:level1deltay} given in terms of $R_1$, $R_2$, and $k$ by
\[
\alpha = \displaystyle\frac{kR_1R_2}{2R_2+R_1},\hspace{.2in} \beta=\frac{k(R_2)^2+(R_2)^2}{2R_2+R_1},\hspace{.1in}\text{and}\hspace{.1in} \gamma=\frac{R_1R_2}{2R_2+R_1}.
\]

\begin{figure}[t]
\setlength{\unitlength}{5cm}
\begin{picture}(1,1)
\put(0,0){\line(3,4){.63}}
\put(1,0){\line(-3,4){.63}}
\put(0,0){\line(1,0){1}}
\put(.375,.84){\line(1,0){.25}}
\put(.45, .85){$R_1$}
\put(0.03,.37){$kR_2$}
\put(.22,.52){$kR_1$}
\put(.6, .74){$R_2$}
\put(.55, .6){$kR_2$}
\put(.2,.03){$R_1$}
\put(.06,.2){$R_2$}
\put(.48,.4){$kR_1$}
\put(.5,0){\line(3,4){.375}}
\put(.5,0){\line(-3,4){.375}}
\put(.125,.5){\line(1,0){.75}}
\end{picture}\hspace{.25in} 
\setlength{\unitlength}{5cm}
\begin{picture}(1,1)
\put(.35,.2){\line(1,0){.3}}
\put(.15,0){\line(1,1){.2}}
\put(.85,0){\line(-1,1){.2}}
\put(.35,.2){\line(0,1){.3}}
\put(.65,.2){\line(0,1){.3}}
\put(.35,.5){\line(-3,1){.1}}
\put(.65,.5){\line(3,1){.1}}
\put(.5,.7){\line(0,-1){.16}}
\put(.5,.54){\line(3,-1){.15}}
\put(.5,.54){\line(-3,-1){.15}}
\put(.5,.7){\line(-3,2){.1}}
\put(.5,.7){\line(3,2){.1}}
\put(.5, .78){$\alpha'$}
\put(.52, .6){$\delta$}
\put(.4,.55){$\epsilon$}
\put(.7,.45){$\zeta$}
\put(.3,.3){$\eta$}
\put(.8,.08){$\theta$}
\put(.5,.15){$\iota$}
\end{picture}\hspace{.5in}\\[6pt]
\setlength{\unitlength}{5cm}
\begin{picture}(1,1)
\put(.3,0){\line(1,2){.2}}
\put(.7,0){\line(-1,2){.2}}
\put(.38,.16){\line(1,0){.24}}
\put(.5,.4){\line(0,1){.2}}
\put(.5,.6){\line(1,1){.15}}
\put(.5,.6){\line(-1,1){.15}}
\put(.55,.7){$\alpha'$}
\put(.52,.5){$\delta$}
\put(.55,.3){$\kappa$}
\put(.5,.1){$\iota$}
\put(.65,.08){$\theta$}
\end{picture}\hspace{.25in} 
\setlength{\unitlength}{5cm}
\begin{picture}(1,1)
\put(.3,0){\line(3,4){.2}}
\put(.7,0){\line(-3,4){.2}}
\put(.5,.25){\line(0,1){.35}}
\put(.5,.6){\line(1,1){.15}}
\put(.5,.6){\line(-1,1){.15}}
\put(.55,.7){$\alpha'$}
\put(.52,.4){$\beta'$}
\put(.65,.1){$\gamma'$}
\end{picture}\hspace{.5in}

\caption{Delta-Y Transformation applied to Level 2 Fractalina}\label{fig:level2deltay}
\end{figure}
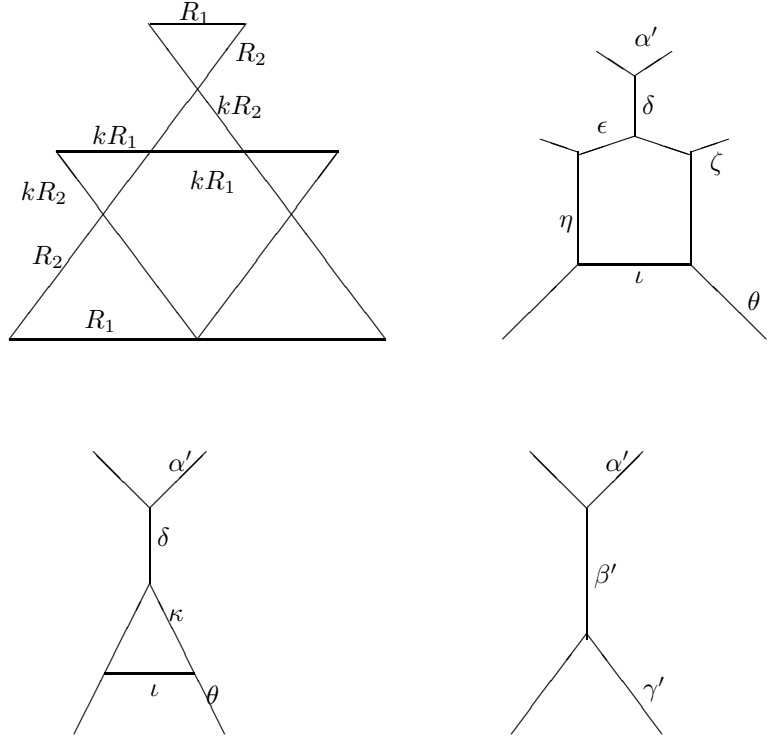

\noindent We now consider the level 2 fractalina.  Figure \ref{fig:level2deltay} shows the unscaled resistances between vertices.  5 applications of a delta-Y transformation yields resistances according to figure \ref{fig:level2deltay}, $\alpha' = \theta = \gamma$, $\delta = \eta = \beta$, $\epsilon = 2\alpha$, and $\iota=2\gamma$. The resistance $\zeta$ is equal to $\alpha$, but this is unimportant as we shall not be considering the effective resistance for any pair of points that a path between which would require these connections.
%
%
A restructuring of the resistances using Kirchhoff's laws gives the lower left side of figure \ref{fig:level2deltay} with the new resistance
\[
\kappa = 2\alpha + \beta =\frac{2kR_1R_2+kR_2^2+(R_2)^2}{2R_2+R_1}.
\]
One more $\Delta$-Y transform results in the resistances
\begin{align*}
 &\beta'=\\ &\frac{(2kR_1R_2+k(R_2)^2+(R_2)^2)^2+(k(R_2)^2+(R_2)^2)(4kR_1R_2+2k(R_2)^2+2(R_2)^2+2R_1R_1)}{(2R_2+R_1)(4kR_1R_2+2k(R_2)^2+2R_2^2+2R_1R_2)}\\[6pt]
&\gamma'=\\&\frac{(2kR_1R_2+k(R_2)^2+(R_2)^2)(2R_2R_1)+(R_1R_2)(4kR_1R_2+2k(R_2)^2+2(R_2)^2+2R_1R_2)}{(2R_2+R_1)(4kR_1R_2+2k(R_2)^2+2(R_2)^2+2R_1R_2)}
\end{align*}


By comparing the effective resistance between $P_1$ and $P_2$ on both levels, we get that $k = 2\gamma/2\gamma'$. Comparing the resistances between $F_3(P_1)$ and $P_1$, we get that $k = (\alpha+\beta+\gamma)/(\alpha'+\beta'+\gamma')$. So now there are two equations and three unknowns: $k$, $R_1$ and $R_2$.  We can reduce the number of unknowns by assuming, without loss of generality,  $R_2=1$. We can do this because resistance is a quadratic form, and we are only interested in the relationship between $R_1$ and $R_2$ up to a constant multiple. This is manifest in the fact that all the above resistances are homogeneous rational functions in $R_1$ and $R_2$. We reduce the first equation

\begin{align}\label{eqn:linares}
\displaystyle\ k=\frac{2kR_1+k+1+R_1}{4kR_1+2k+R_1}
\end{align}

\noindent Solve for $R_1$ in terms of k,
\[
R_1=\frac{1-k-2k^2}{4k^2-k-1}
\]
Assume $k=(3+\sqrt{41})/16$.   This implies that $R_1=(-1+\sqrt{41})/4$. Of course, it must be verified that these values satisfy the second equation, which it also satisfies. Thus the scaling factor of the resistance form on the fractalina fractal is $\displaystyle\frac{3+\sqrt{41}}{16}$.
%
\section{Pillow Fractal}\label{sec:pillow}
 
\begin{figure}[t]
\setlength{\unitlength}{3.5cm}
\begin{picture}(1,1)
\put(0,0){\line(0,1){.75}}
\put(0,0){\line(1,0){.75}}
\put(0,0){\line(1,0){1}}
\put(1,0){\line(0,1){.75}}
\put(0,.75){\line(1,0){1}}
\put(.35,-.15){$\bf F_0$}
\put(0,0){\circle*{.04}}
\put(1,.75){\circle*{.04}}
\put(1,0){\circle*{.04}}
\put(0,.75){\circle*{.04}}
\end{picture} \hspace{.5in}
\setlength{\unitlength}{3.5cm}
\begin{picture}(1,1)
\put(0,0){\line(0,1){.75}}
\put(0,0){\line(1,0){.5}}
\put(.5,0){\line(0,1){.75}}
\put(0,.75){\line(1,0){.5}}
\put(0,0){\circle*{.04}}
\put(.5,.75){\circle*{.04}}
\put(.5,0){\circle*{.04}}
\put(0,.75){\circle*{.04}}
\put(.7,0){$\bf\leftrightarrow$}
\put(.7, .75){$\leftrightarrow$}
\end{picture} 
\setlength{\unitlength}{3.5cm}
\begin{picture}(1,1)
\put(0,0){\line(0,1){.75}}
\put(0,0){\line(1,0){.5}}
\put(.5,0){\line(0,1){.75}}
\put(0,.75){\line(1,0){.5}}
\put(0,0){\circle*{.04}}
\put(.5,.75){\circle*{.04}}
\put(.5,0){\circle*{.04}}
\put(0,.75){\circle*{.04}}
\end{picture}

\begin{picture}(1,1)
\put(0,0){\line(0,1){.75}}
\put(0,0){\line(1,0){1}}
\put(1,0){\line(0,1){.75}}
\put(0,.75){\line(1,0){1}}
\put(0,0){\circle*{.04}}
\put(0, .75){\circle*{.04}}
\put(1, 0){\circle*{.04}}
\put(.5,0){\circle*{.04}}
\put(.5,.75){\circle*{.04}}
\put(1, .75){\circle*{.04}}
\put(.5,0){\qbezier(0.0,0.0)(.1,.37)
(0,.75)}
\put(.5,0){\qbezier(0.0,0.0)(-.1,.37)
(0,.75)}
\put(.5,-.15){$\bf F_1$}
\end{picture}
\caption{Contrustion of $F_1$ from $F_0$ }\label{fig:construction}
\end{figure}
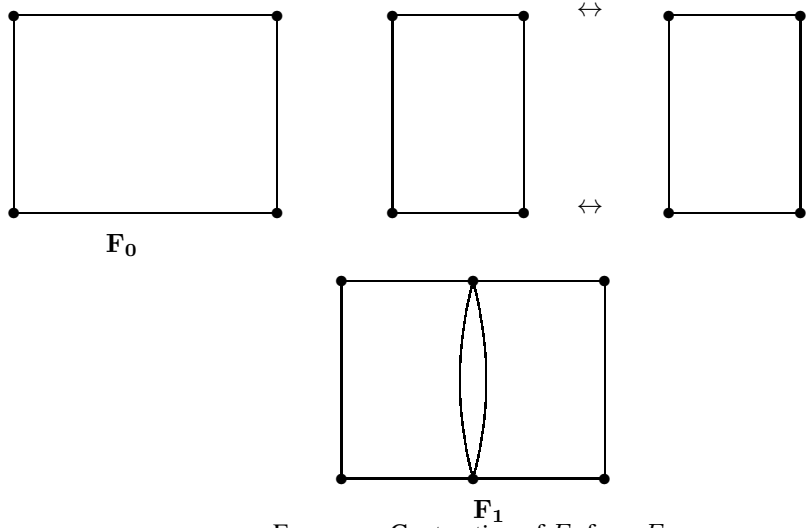

The pillow fractal is a self-affine fractal in that the contraction mappings that produce this fractal are not rigid. The basic building blocks of the pillow are rectangles and there are two contraction maps. Both contraction maps rotate the fractal 90 degrees, gluing together two of the corners while leaving the edges separate. Figure \ref{fig:construction} shows this construction.

In \cite{NT08}, the pillow is constructed as the limit space of a self-similar group. A recurrence relation is determined from the action of this group, which gives equations to determine the resistance of the Schreier graphs approximating the limit space. These Schreier grpahs are different than the graphs used used in the current paper to approximate the pillow, however, the results are consistent.

The energy on the pillow is built up from approximating graphs, the 0th level is a ``rectangle'' in that it has four nodes $F_0 =\set{v_1,v_2,v_3,v_4}$, a node at every corner.  We are assuming rectangular symmetry of the pillow. In particular the energy form is invariant under horizontal and vertical flips as well as 180 degree rotation. More precisely, for the maps $\sigma = (v_1v_2)(v_3v_4)$ , $\tau = (v_1v_3)(v_2v_4)$, $\rho = (v_1v_4)(v_2v_3)$ in the symmetric group of $F_0$,
\[
\eng_0(f) = \eng_0(f\circ\sigma) = \eng_0(f\circ\tau) = \eng_0(f\circ \rho)
\] 
This invariance implies that we can assume that $C_{v_1v_3} = C_{v_2v_4} := C_1$, $C_{v_1v_2} = C_{v_3v_4} := C_2$, and $C_{v_1v_4} = C_{v_2v_3} := C_3$.

Define $F_1 = F_0\cup\set{u_1,u_2}$, and the first contraction maps $v_1 \to v_1$, $v_2 \to u_2$, $v_3\to v_2$, and $v_4\to u_2$, while the second maps $v_1 \to v_1$, $v_2 \to u_2$, $v_3\to v_2$ and $v_4\to u_2$. The conductances of $F_1$ are determined by these contractions, noting that we get two edges connecting $u_1$ to $u_2$ with conductance $C_1$ (equivalently, we can think of $C_{u_1u_2} = 2C_1$). Figure \ref{fig:pillowconductance} shows $F_0$ and $F_1$.

The rectangular symmetry from $F_0$ can be extended to $F_1$, that is that $\eng_2$ is invariant under $\tilde{\sigma} = (v_1v_2)(v_3v_4)(u_1u_2)$, $\tilde{\tau} = (v_1v_3)(v_2v_4)$, $\tilde{\rho} = (v_1v_4)(v_2v_3)(u_1u_2)$. Notice that, for example $\tilde{\sigma}|_{F_0} = \sigma$

\begin{theorem}
The resistance scaling factor between any two levels, $\mathscr{E}_i$ and $\mathscr{E}_{i+1}$ is $\sqrt[3]2$
\end{theorem}

The rest of this section is the proof of this. To find the resistance scaling for the pillow, we employ a different strategy than in section \ref{sec:fractalina}.  We look at specific states of the fractal, and see how they correspond to the graph at each level by assigning potentials to each node.

To compare states consistently between graphs, we take a function $f:F_0\to \Bbb{R}$, we get a function on $F_1$ by \emph{harmonic extension} --- the function $\tilde{f}$ such that $\tilde{f}|_{F_0} = f$, and
\begin{align}
\label{eq:harmonic}\tilde{f}(p) = \sum_{q\in F_1}C_{pq}\tilde{f}(q)\Big/\sum_{q\in F_1}C_{pq}
\end{align}
for all verteces $p\notin F_0$, and $C_{pq}$ is the conductance between $p$ and $q$. The harmonic extension $\tilde{f}$ is the energy minimizing extension of $f$. That is to say $\eng_1(\tilde{f}) \leq \eng_1(g)$ for any state $g$ such that $g|_{F_0} = f$, or that $\eng_1(\tilde{f}) = \trace \eng_1(f)$. For more on this see \cite{Str06,DS84}. 

We make some important observations about harmonic extension that are straight forward to verify. First, harmonic extension is linear: if $\tilde{f}$ and $\tilde{g}$ are the harmonic extensions of $f$ and $g$ respectively, then $\tilde{f}+c\tilde{g}$ if the harmonic extension of $f+cg$, for any constant $c$. Second, the extension of a constant function is constant. This is easy to see, because the energy of a constant function is 0, so the constant function is an energy minimizing extension. Finally, harmonic extension is invariant under symmetries of $F_0$ and $F_1$, that is to say if $f\circ \sigma = f$, then $\tilde{f}\circ\tilde{\sigma} = \tilde{f}$ (similarly for $\tau$ and $\rho$), this can be seen by examining the definition in \ref{eq:harmonic}.

Comparing the energy of the state $f$ to the energy of $\tilde{f}$ will give the resistance scaling factor, $\eng_0(f)/\eng_1(\tilde{f})$. However, this calculation leaves unknown quantities. We will have to consider three functions $f_1,f_2,f_3:F_0\to \Bbb{R}$ to solve for these quantities, and obtain the scaling factor.

\begin{figure}[t]
\setlength{\unitlength}{4.5cm}
\begin{picture}(1,1)
\put(0,0){\line(0,1){.75}}
\put(0,0){\line(4,3){1}}
\put(1,0){\line(-4,3){1}}
\put(0,0){\line(1,0){.75}}
\put(0,0){\line(1,0){1}}
\put(1,0){\line(0,1){.75}}
\put(0,.75){\line(1,0){1}}
\put(-.07,-.07){$v_1$}
\put(-.07, .78){$v_2$}
\put(1, -.07){$v_3$}
\put(1, .78){$v_4$}
\put(.5, .78){$C_1$}
\put(.5, -.07){$C_1$}
\put(1.01, .5){$C_2$}
\put(-.1,.5){$C_2$}
\put(.75, .5){$C_3$}
\put(.35,-.15){$ F_0$}
\end{picture} \hspace{.5in}
\setlength{\unitlength}{4.5cm}
\begin{picture}(1,1)
\put(0,0){\line(2,3){.5}}
\put(0.5,0){\line(2,3){.5}}
\put(0.5,0){\line(-2,3){.5}}
\put(1,0){\line(-2,3){.5}}
\put(0,0){\line(0,1){.75}}
\put(0,0){\line(1,0){1}}
\put(1,0){\line(0,1){.75}}
\put(0,.75){\line(1,0){1}}
\put(-.07,-.07){$v_1$}
\put(-.07, .78){$v_2$}
\put(1, -.07){$v_3$}
\put(1, .78){$v_4$}
\put(.5, .78){$u_1$}
\put(.5, -.07){$u_2$}

\put(.5,0){\qbezier(0.0,0.0)(.1,.37)
(0,.75)}
\put(.5,0){\qbezier(0.0,0.0)(-.1,.37)
(0,.75)}
\put(.46,.35){$C_1$}
\put(.75, .78){$C_2$}
\put(.75, -.07){$C_2$}
\put(1.01, .5){$C_1$}
\put(-.1,.5){$C_1$}
\put(.75, .5){$C_3$}
\put(.5,-.15){$ F_1$}
\end{picture}\label{fig:Conductances}
\caption{$F_0$ and $F_1$ Conductances}\label{fig:pillowconductance}
\end{figure}

\emph{Configuration 1} Define $f_1(v_1) = f_1(v_2) = f_1(v_3) = 0$ and $f_1(v_4)=1$

To ensure that we get the appropriate energy we calculate the energy minimizing potentials at the unaccounted for verteces of $F_1$, $\tilde{f}_1(u_1) = x$ and $\tilde{f}_1(u_2) = y$ as in figure \ref{fig:config1}. The values $x$ and $y$ satisfy equation \ref{eq:harmonic}. Written in terms of $C_i$, for $i=1,2,3$,
\begin{align}
x=\dfrac{2C_1y+C_2}{2(C_1+C_2+C_3)}\quad\text{and}\quad y=\dfrac{2C_1x+C_3}{2(C_1+C_2+C_3)}.\label{xy}
\end{align}
These equations are used to calculate the energy form for the first and second orders, whose ratio is
\[
\dfrac{\mathscr{E}_0(f_1)}{\mathscr{E}_1(\tilde{f}_1)}=\dfrac{C_1+C_2+C_3}{C_1(1+2(x-y)^2)+C_2(x^2+2y^2+(x-1)^2)+C_3(2x^2+y^2+(y-1)^2)}.
\]

We complete similar calculations for two more cases.\\

\noindent\emph{Configuration 2}: Define $f_2(v_2)=f_2(v_3)=0$ and $f_2(v_1)=f_2(v_4)=1$. Notice that $f_2\circ \rho = f_2$, so $\tilde{f_2}(u_2) = \tilde{f}_2(\tilde{\rho}(u_1)) = \tilde{f}_2(u_1)$. Furthermore, $f_2 + f_2\circ \tau$ is the constant function 1, so $\tilde{f}_2(u_1)+\tilde{f}_2(\tilde{\tau}(u_1)) = 2\tilde{f}_2(u_1) = 1$. Thus $\tilde{f}(u_1) = \tilde{f}_2(u_2) = 1/2$.

Since we know all the values of $\tilde{f}_2$, we calculate
\[
\dfrac{\mathscr{E}_0(f_2)}{\mathscr{E}_1(\tilde{f}_2)}=\dfrac{2C_1+2C_2}{2C_1+C_2+C_3}
\]
\noindent\emph{Configuration 3}: Define $f_3(v_1)=f_3(v_2)=0$ and $f_3(v_3)=f_3(v_4)=1$. Similar to the computation of $\tilde{f}_2$, $f_3\circ \sigma = f_3$ and $f_3 + f_3\circ \tau \equiv 1$, so $\tilde{f}_3(u_1) = \tilde{f_3}(u_2) = 1/2$. Leaving us to compute
\[
\dfrac{\mathscr{E}_0(f_3)}{\mathscr{E}_1(\tilde{f_3})}=\dfrac{2C_1+2C_3}{C_2+C_3}
\]

\begin{figure}
\setlength{\unitlength}{5cm}
\begin{picture}(1,1)
\put(0,0){\line(0,1){.75}}
\put(0,0){\line(4,3){1}}
\put(1,0){\line(-4,3){1}}
\put(0,0){\line(1,0){.75}}
\put(0,0){\line(1,0){1}}
\put(1,0){\line(0,1){.75}}
\put(0,.75){\line(1,0){1}}
\put(-.07,-.07){$0$}
\put(-.07, .78){$0$}
\put(1, -.07){$0$}
\put(1, .78){$1$}
\put(.5, .78){$C_1$}
\put(1.01, .5){$C_2$}
\put(.75, .5){$C_3$}
\end{picture}\hspace{.5in}
\setlength{\unitlength}{5cm}
\begin{picture}(1,1)
\put(0,0){\line(2,3){.5}}
\put(0.5,0){\line(2,3){.5}}
\put(0.5,0){\line(-2,3){.5}}
\put(1,0){\line(-2,3){.5}}
\put(0,0){\line(0,1){.75}}
\put(0,0){\line(1,0){1}}
\put(1,0){\line(0,1){.75}}
\put(0,.75){\line(1,0){1}}
\put(-.07,-.07){$0$}
\put(-.07, .78){$0$}
\put(1, -.07){$0$}
\put(1, .78){$1$}
\put(.5,0){\qbezier(0.0,0.0)(.1,.37)
(0,.75)}
\put(.5,0){\qbezier(0.0,0.0)(-.1,.37)
(0,.75)}
\put(.46,.35){$C_1$}
\put(.75, .78){$C_2$}
\put(.75, -.07){$C_2$}
\put(1.01, .5){$C_1$}
\put(-.1,.5){$C_1$}
\put(.75, .5){$C_3$}
\put(.5, .78){$x$}
\put(.5, -.07){$y$}
\end{picture}
\caption{Configuration 1}\label{fig:config1}
\end{figure}
Assume $\dfrac{\mathscr{E}_0}{\mathscr{E}_1}=2^{1/3}$ for the latter two equations, and without loss of generality, let $C_1 = 1$. This implies that
\[
2^{1/3}=\dfrac{2+2C_2}{2+C_2+C_3}\quad\text{and}\quad2^{1/3}=\dfrac{2+2C_3}{C_2+C_3}
\]
Solve for the two conductances we get that $C_3=2^{-1/3}$ and $C_2=2^{-1/3}+2^{1/3}$. It now must be verified that these values satisfy the equation for $\mathscr{E}_0(f_1) / \mathscr{E}_1(\tilde{f}_1)$. Apply the values of $C_i$ to the equations \ref{xy}, we get that $x = 2^{-5/3}$ and $y = 2^{-1}-2^{-5/3}$. Applying these results into the equation above we get that $\dfrac{\mathscr{E}_0(f_1)}{\mathscr{E}_1(\tilde{f}_1)} = 2^{-1/3}$, which proves consistency.
%
\bibliography{Fractals}{}
\bibliographystyle{amsalpha}
%
\end{document}